\documentclass{IEEEtran4PSCC}
\ifCLASSINFOpdf
\else
\fi
\usepackage[cmex10]{amsmath}
\usepackage{xspace,framed,xcolor}

\usepackage{amssymb}
\usepackage{amsmath}
\usepackage{amsthm}
\usepackage[english]{babel}
\usepackage{xspace}
\usepackage{cases}
\usepackage{color}
\usepackage{graphicx}
\usepackage[caption=false]{subfig}
\usepackage{epstopdf}
\usepackage{upgreek}
\usepackage{enumitem}
\usepackage{booktabs}
\usepackage{array}
\usepackage{blindtext}
\usepackage{setspace}
\usepackage[hidelinks,bookmarks=false]{hyperref}

\newtheorem{thm}{Theorem}

\newcommand{\nbus}{n}

\newcommand{\setbus}{\mathcal{N}}

\newcommand{\nline}{m}

\newcommand{\basisfun}{\psi}

\newcommand{\adjustlength}{-4.1mm}

\renewcommand{\arraystretch}{1.15}
\newcommand{\pdf}{\textsc{pdf}\xspace}
\newcommand{\pdfs}{\textsc{pdf}s\xspace}

\newcommand{\opf}{\textsc{opf}\xspace}

\newcommand{\sopf}{\ccopf}
\newcommand{\dc}{\textsc{dc}\xspace}
\newcommand{\ac}{\textsc{ac}\xspace}

\newcommand{\dcopf}{\textsc{dc}-\textsc{opf}\xspace}

\newcommand{\ccopf}{cc\textsc{opf}\xspace}
\newcommand{\hopf}{h\textsc{opf}\xspace}

\newcommand{\pce}{\textsc{pce}\xspace}

\newcommand{\qp}{\textsc{qp}\xspace}

\newcommand{\bs}[1]{\boldsymbol{#1}}
\newcommand{\rv}[1]{\mathsf{#1}}

\newcommand{\ev}[1]{\mathbb{E}\left[#1\right]}

\newcommand{\pr}[1]{ {{\mathbb{P}}}[#1] }

\newcommand{\pvar}[1]{p^{g}_{#1}}
\newcommand{\qvar}[1]{q^{g}_{#1}}
\newcommand{\pdvar}[1]{p^{d}_{#1}}
\newcommand{\qdvar}[1]{q^{d}_{#1}}
\newcommand{\rvpvar}[1]{\rv{p}^{g}_{#1}}
\newcommand{\rvqvar}[1]{\rv{q}^{g}_{#1}}
\newcommand{\rvpdvar}[1]{\rv{p}^{d}_{#1}}
\newcommand{\rvqdvar}[1]{\rv{q}^{d}_{#1}}

\newcommand{\kkt}{\textsc{kkt}\xspace}

\newcommand{\caseunc}{\textsc{c}1\xspace}
\newcommand{\casecon}{\textsc{c}2\xspace}

\newlist{mylist}{enumerate}{1}
\setlist[mylist]{label=\textsc{a}\oldstylenums{\arabic*}}
\newlist{mylistlist}{enumerate}{1}
\setlist[mylistlist]{label=\textsc{q}\oldstylenums{\arabic*}}

\definecolor{till}{rgb}{.1,.4,.9}
\definecolor{timm}{rgb}{.0,.0,.0} 
\definecolor{veit}{rgb}{.8,.1,.1}

\newlist{mypcelist}{enumerate}{1}
\setlist[mypcelist]{label=\textsc{p}\oldstylenums{\arabic*}}
\newlist{myinglist}{enumerate}{1}
\setlist[myinglist]{label=\textsc{i}\oldstylenums{\arabic*}}
\usepackage{pgf, algorithm, algpseudocode}

\usepackage[citestyle=numeric-comp, uniquename=false, sorting=none, maxbibnames=6, maxcitenames=2, doi=false, url=false, isbn=false, eprint=false, backend=bibtex]{biblatex}
\begin{document}
%
\title{The Price of Uncertainty:\\ Chance-constrained \opf vs. In-hindsight \opf}

\author{
\IEEEauthorblockN{Tillmann M\"uhlpfordt,  Veit Hagenmeyer, Timm Faulwasser}
\IEEEauthorblockA{Institute for Automation and Applied Informatics, Karlsruhe Institute of Technology\\ 
Karlsruhe, Germany\\
\{tillmann.muehlpfordt, veit.hagenmeyer, timm.faulwasser\}@kit.edu}
}


\maketitle

\begin{abstract}
The operation of power systems has become more challenging due to feed-in of volatile renewable energy sources.
Chance-constrained optimal power flow (\ccopf) is one possibility to explicitly consider volatility via probabilistic uncertainties resulting in mean-optimal feedback policies.
These policies are computed before knowledge of the realization of the uncertainty is available.
On the other hand, the hypothetical case of computing the power injections knowing every realization beforehand---called in-hindsight \opf (\hopf)---cannot be outperformed w.r.t. costs and constraint satisfaction.
In this paper, we  investigate how \ccopf feedback relates to the full-information \hopf.
To this end, we introduce different dimensions of the price of uncertainty.
Using mild assumptions on the uncertainty we present sufficient conditions when \ccopf is identical to \hopf.
We suggest using the total variational distance of probability densities to quantify the performance gap of \hopf and \ccopf.
Finally, we draw upon a tutorial example to illustrate our results. 
\end{abstract}


\thanksto{The authors acknowledge support  by the  Helmholtz Association under the Joint Initiative ``Energy System 2050 -- A Contribution of the Research Field Energy''.}

\section{Introduction}
Power generation currently undergoes a paradigm shift:
The tremendous installation of renewable energy generation units leads to varying operating conditions, fluctuating reserve capacities, and increasingly changing power flows on different levels of the grid, among others.
As the investment in sustainable energy production continues, operational control strategies have to be developed allowing secure grid operation despite the presence of uncertainties.
Worst-case methods \cite{Warrington13} and stochastic programming methods \cite{Bouffard08} have been proposed to deal with uncertainties in power systems applications on the secondary and tertiary control level.
The present paper focuses on chance-constrained optimal power flow (\sopf) formulations \cite{Zhang10,Vrakopoulou12,Bienstock14,Roald13}.
These formulations offer a framework to compute power injections that guarantee power balance, while respecting generation and transmission limits with a user-specified probability.
Much of the literature on \sopf focuses on the \dc setting \cite{Vrakopoulou12,Bienstock14,Roald13,Roald15b,Muehlpfordt17a}.
Recently, however, the \ac case has been considered too \cite{Roald17,Muehlpfordt16b,Vrakopoulou13}.

The setup of \sopf is conceptually appealing, because it extends conventional automatic generation control to the uncertain case; control policies rather than single generator set points are determined.
The parameters of the control policies are the decision variables of \sopf problems.
Immediately upon measuring a specific realization, the optimal control policies  from \sopf provide power injections which satisfy the generation and transmission limits in the usual chance constraint sense.
However, chance-constrained optimization problems can be intrinsically difficult to solve.
This is why much effort has been put into reformulating chance constraints as deterministic (convex) constraints that lead to tractable optimization problems \cite{Bienstock14,Roald15b,Vrakopoulou13c}.
Effectively, chance constraint reformulations introduce constraint tightening; the so-called \emph{uncertainty margin} \cite{Qu15}.
Constraint tightening, in turn, may lead to conservative \sopf solutions, which may then result in higher (expected) operational costs of the power system.
In other words, as claimed in  \cite[p. 5]{Kall94}: ``The stochastic solution [\ldots] is normally never optimal [\ldots], at the same time, it is also hardly ever really bad.''

Chance-constrained optimal power flow is fundamentally challenging because control policies have to be determined \emph{before} the realization of the uncertainty is known.
This induces a henceforth called ``price of uncertainty'' that has---as the above considerations show---at least four dimensions in the context of power flow problems:\footnote{Relative priorities of these dimensions are mostly case-specific and thus in general hard to state.}
\begin{itemize}
	\item \emph{Cost}: What additional (expected) monetary costs are induced due to uncertainties?
	\item \emph{Operation}: How does the grid operation change in the presence of uncertainties?	
	\item \emph{Computation}: Is the optimization problem from \ccopf computationally tractable?
	\item \emph{Feasibility}: Are the operating limits satisfied for all realizations?
\end{itemize}
In contrast to \sopf that considers stochasticity of the uncertainty explicitly via random variables, online optimization approaches to power flow problems \cite{Hauswirth17,DallAnese16,Gan16} react directly to the \emph{actual realization} of the uncertainty.
	Online optimization approaches hence admit an intuitive interpretation as (implicitly defined) feedback controllers acting on the physical grid, which becomes the controlled system \cite{Hauswirth17,Gan16}.
	We remark that in the last 25 years the systems and control community has witnessed vast research efforts in the field of feedback control via numerical online optimization, i.e. model predictive control; moreover, it has changed industrial control practice \cite[p. xi]{maciejowski2002predictive}.
	When online optimization is applied to power systems, it consequently brings about the challenges known from model predictive control: state estimation, and immediate, accurate and reliable solution of (non-convex) optimization problems, and implementation thereof on existing automation systems.

In view of the above, it is fair to ask for the potential losses/gains of control policies obtained via \sopf compared to  online optimization.
To the end of providing first elements of an answer, herein we introduce fully-informed in-hindsight \opf (\hopf).
This refers to the unrealistic case that the solution of the \opf problem is immediately known \emph{for all} realizations of the uncertainty, whereby
each individual \hopf solution is per-sample optimal.
Per-sample optimality here means that for the respective realization the minimum cost is attained, and the constraints are strictly satisfied.
In contrast to \sopf, the constraints are never violated with \hopf; to quote again \cite[p. 5]{Kall94}: ``The \textsc{iq} of hindsight is very high.''
Note that \hopf contains the solution via online optimization.

The contribution of the present paper is to answer the question when \sopf and \hopf provide equivalent solutions.
In that case, the reaction to any realization of the uncertainty is computed in a single numerical run \emph{and} known to be per-sample optimal.
An online optimization is then unnecessary and can be replaced by a single control policy that is computed once offline.
In case of \sopf and \hopf providing different solutions, we suggest using the total variational distance as a metric to quantify the price of uncertainty.
Our findings are illustrated by means of a tutorial three-bus system.

\section{Problem Setup}
Consider a power system with $\nbus$ buses and $\nline$ lines.
The bus index set is $\setbus = \{1, \hdots, n\}$ with $| \setbus | = \nbus$.
For ease of presentation, we assume that a single generator and a single load are connected to every bus.
Additionally, the load is uncertain in the sense that
the actual value of the load is unknown.
Hence, the load is modeled as a continuous second-order random variable, i.e. $\rvpdvar{i} \in \mathrm{L}^2(\Omega_i, \mu_i; \mathbb{R})$ for every bus $i \in \setbus$.\footnote{%
The space $\mathrm{L}^2(\Omega_i, \mu_i; \mathbb{R})$ is the Hilbert space of second-order $\mathbb{R}$-valued random variables with support $\Omega_i$ and probability measure $\mu_i$ \cite{Sullivan15book}.
For simplicity, we assume the support of $\mu_i$ is equal to the sample space $\Omega_i$.}
To ensure power balance in the presence of uncertainties, at least one controllable active power injection must then necessarily be modeled as a random variable, too.
In other words, if a load changes unexpectedly, some generator has to change its injection too.
To simplify presentation, we assume that all controllable generators might in principle react to the uncertainties; thus, every generator is modeled as a random variable.\footnote{\label{footnote:TrivialUncertainty}If a generator/load is \emph{not} uncertain, it can be modeled as a random variable with a Dirac-delta probability density centered around the deterministic value.}
Then, the net active and net reactive power for bus $i$ are modeled as the following random variables
\begin{subequations}
	\label{eq:NetPower_rv}	
	\begin{align}
	\rv{p}_i &= \rvpvar{i} + \rvpdvar{i} \in \mathrm{L}^2(\Omega_i, \mu_i; \mathbb{R}), \\
	\rv{q}_i &= \rvqvar{i} + \rvqdvar{i} \in \mathrm{L}^2(\Omega_i, \mu_i; \mathbb{R}).
	\end{align}
\end{subequations}
For each realization of the random variables in \eqref{eq:NetPower_rv},
\begin{subequations}
	\label{eq:NetPower_realization}
	\begin{align}
	p_i &= \pvar{i} + \pdvar{i} \in \mathbb{R}, \\
	q_i &= \qvar{i} + \qdvar{i} \in \mathbb{R},
	\end{align}
\end{subequations}
the optimal power flow (\opf) problem can be formulated as
\begin{subequations}
	\label{eq:AC_OPF}
	\begin{align}
	\underset{\pvar{}, \qvar{}}{\min}\quad  & \sum_{i \in \setbus} c_i(\pvar{i})  \\
	\mathrm{s.\, t.} \, ~ ~ ~ \,
	\label{eq:PFE}
	& g(\pvar{},\qvar{},v,\theta; \pdvar{}, \qdvar{}) = 0,\\
	\label{eq:CON_PV}
	& x_i^{\text{min}} \leq x_i \leq x_i^{\text{max}},  && \forall i \in \setbus,\\
	& \nonumber \forall x_i \in \{ \pvar{i}, \, \qvar{i}, v_i, \theta_i  \}, \\
	& \theta_{i_s} = 0, && i_s  \in \setbus,\\
	\label{eq:CON_i}
	& i_{i,j}^{\text{min}} \leq i_{i,j} \leq i_{i,j}^{\text{max}}, && \forall i,j \in \setbus,
	\end{align}
\end{subequations}
where $v_i$, $\theta_i$ are the magnitude, phase of the voltage phasor at bus~$i$, respectively, and $i_{i,j}$ is the magnitude of the transmitted current between buses~$i$ and~$j$.
Problem~\eqref{eq:AC_OPF} minimizes the sum of active power generation costs $c_i$ subject to the power flow equations \eqref{eq:PFE}, generation constraints and voltage constraints \eqref{eq:CON_PV}, the slack constraint, and transmission constraints \eqref{eq:CON_i}.
For simplicity, the high-voltage solution to Problem \eqref{eq:AC_OPF} is assumed to exist for all realizations of \eqref{eq:NetPower_rv}.
The minimizer of the \opf Problem \eqref{eq:AC_OPF} depends on the specific realization of the power demands $\pdvar{}$ and $\qdvar{}$
\begin{equation}
\label{eq:argmin_operator}
\begin{bmatrix}
\pvar{}\phantom{}^\star \\
\qvar{}\phantom{}^\star
\end{bmatrix}
=
\begin{bmatrix}
\nu_p(\pdvar{}, \qdvar{}) \\
\nu_q(\pdvar{}, \qdvar{}) 
\end{bmatrix}
 := \underset{\pvar{}, \qvar{}}{\operatorname{argmin}} \: \, \text{Problem \eqref{eq:AC_OPF}},
\end{equation}
where $\nu_p, \nu_q: \mathbb{R}^{2 \nbus} \rightarrow \mathbb{R}^{\nbus}$ comprise the argmin operator of the \opf Problem \eqref{eq:AC_OPF}.
In-hindsight \opf refers to the unrealistic situation, in which the solution to the \opf Problem~\eqref{eq:AC_OPF} is known and immediately available \emph{for all} realizations of \eqref{eq:NetPower_rv}.
Consequently, the result of \hopf is itself a random variable.
More precisely, it is the random-variable optimal active and reactive power generation $\rv{p}^{g \star}, \rv{q}^{g \star}$.
Algorithm~\ref{fig:Algorithm_hopf} summarizes the pseudo code of \hopf.
\begin{algorithm}
	\setstretch{1.2}
	\vspace{2mm}
	\centering
	\small
	\begin{algorithmic}[1]
		\State Choose number of samples $N \in \mathbb{N}$.
		\State Draw $N$ samples $\{(\pdvar{})_k, (\qdvar{})_k \}_{k=1}^N$ from $\rvpdvar{}, \rvqdvar{}$.
		\For{$k = 1, \hdots, N$}
		\State Pick $k$\textsuperscript{th} sample $(\pdvar{})_k, (\qdvar{})_k$.
		\State Solve 
		$$
		\begin{bmatrix}
		(\pvar{}\phantom{}^\star)_k\\
		(\qvar{}\phantom{}^\star)_k
		\end{bmatrix}
		=
		\begin{bmatrix}
		\nu_p((\pdvar{})_k, (\qdvar{})_k) \\
		\nu_q((\pdvar{})_k, (\qdvar{})_k) 
		\end{bmatrix}.
		$$
		\EndFor
		\State Result: $\rv{p}^{g \star}, \rv{q}^{g \star}$.
	\end{algorithmic}
	\normalsize
	\caption{Description of \hopf.}
	\label{fig:Algorithm_hopf}
\end{algorithm}

In-hindsight \opf is an important limiting case with the following properties:
\begin{enumerate}
	\item \label{item:OnlineOptimization} The solution of \hopf includes the optimal realization $p^{g \star}$, $q^{g \star}$ of power injections for the a priori unknown \emph{actual} realization $\pdvar{}, \qdvar{}$ of random generation/demand from \eqref{eq:NetPower_rv}.
	\item Every realization of the \hopf solution $\rv{p}^{g \star}, \rv{q}^{g \star}$ satisfies the power flow equations and the inequality constraints.
	\item \hopf provides the best distribution of optimal costs in the sense that every sample solution is known to be optimal.
	\item \hopf provides the best distribution of optimal active power generations that are always feasible.
\end{enumerate}
Note that item \ref{item:OnlineOptimization}) corresponds to the situation from online optimization approaches to \opf problems \cite{Hauswirth17,DallAnese16,Gan16}: assuming the state of the grid is immediately and accurately available, the \opf Problem \eqref{eq:AC_OPF} is solved online to provide power injections that balance the grid and satisfy generation and transmission limits.\footnote{As mentioned in the introduction, the idea is similar to model predictive control, except that single-stage \opf problems do not involve any dynamics.}
This way, stochasticity of uncertainties does not have to be accounted for explicitly, because the online optimization algorithm reacts to any of its effects, ideally, in real time.

In-hindsight \opf has conceptual difficulties: The naive algorithm described in Algorithm~\ref{fig:Algorithm_hopf} results in a mere look-up table, i.e. how the optimal solutions depends on the respective realization is not immediate by means of a function/feedback law.\footnote{In case of \dcopf, multiparametric programming techniques \cite{Bemporad02a,Vrakopoulou17} overcome this issue;
yet, this might lead to implementation difficulties. A thorough investigation is beyond the scope of this paper.}
	In case one resorts to fast online solutions based on measured disturbances, the control policy is defined implicitly via optimization, which complicates its analysis.
	The online optimization further has to be implemented using existing control hardware, which poses another challenge.

This motivates a different approach, namely aforementioned chance-constrained optimal power flow (\sopf).
It alleviates the conceptual disadvantages of \hopf mentioned above: in a single numerical run the generation response to all load fluctuations is obtained.
This is achieved by optimizing over control policies.
The \sopf problem can be formulated as follows \cite{Roald13,Muehlpfordt17a,Roald17}
\begin{subequations}
	\label{eq:AC_sOPF}
	\begin{align}
	\label{eq:AC_sOPF_cost}
	\underset{\alpha_p, \alpha_q}{\min}\quad  & \mathbb{E}\Big[\sum_{i \in \setbus} c_i(\rvpvar{i}) \Big] \\
	\mathrm{s.\, t.} \, ~ ~ ~ \,
	\label{eq:PFE_RV}
	& g(\rvpvar{},\rvqvar{},\rv{v},\uptheta; \rvpdvar{}, \rvqdvar{}) = 0,\\
	\label{eq:CON_CC1}
	& \pr{x_i^{\text{min}} \leq \rv{x}_i} \geq 1 - \varepsilon, && \forall i \in \setbus,\\
	\label{eq:CON_CC2}
	& \pr{ \rv{x}_i \leq x_i^{\text{max}}} \geq 1 - \varepsilon, && \forall i \in \setbus,\\
	\label{eq:CON_CC3}
	& \nonumber\forall \rv{x}_i \in \{ \rvpvar{i},\, \rvqvar{i},\, \rv{v}_i,\, \uptheta_i \}, \\
	& \uptheta_{i_s} = 0, && i_s \in \setbus,\\
	\label{eq:CON_CC_thermal1}
	& \pr{ i_{i,j}^{\text{min}} \leq \rv{i}_{i,j}} \geq 1 - \varepsilon, && \forall i, j \in \setbus,\\
	\label{eq:CON_CC_thermal2}
	& \pr{  \rv{i}_{i,j} \leq i_{i,j}^{\text{max}}} \geq 1 - \varepsilon, && \forall i, j \in \setbus,\\
	\label{eq:Policy_p}
	& \rvpvar{} = \eta_p(\rvpdvar{}, \rvqdvar{}; \alpha_p), \\
	\label{eq:Policy_q}
	& \rvqvar{} = \eta_q(\rvpdvar{}, \rvqdvar{}; \alpha_q),
	\end{align}
\end{subequations}
where $\mathbb{E}[\cdot]$ is the expected value, and all sans-serif symbols denote random variables.
We remark that minimizing the expected value of the cost as in \eqref{eq:AC_sOPF_cost} is a modeling choice. We refer to \cite{Muehlpfordt17a} for other choices; for instance, the objective may also include terms penalizing the variance of the cost function.
The formulation \eqref{eq:PFE_RV} of the power flow equations in terms of random variables ensures that the power balance holds for all realizations of the uncertainty \cite{Muehlpfordt16b,Muehlpfordt17a}, assuming power flow feasibility of the high-voltage solution.
The inequality constraints for generation and voltage limits are reformulated in \eqref{eq:CON_CC1}--\eqref{eq:CON_CC3} as single-sided chance constraints, which is again a modeling choice.
Similarly, the transmission limits are modeled as single-sided chance constraints \eqref{eq:CON_CC_thermal1}, \eqref{eq:CON_CC_thermal2}.
The control policies $\eta_p$, $\eta_q$ are introduced in \eqref{eq:Policy_p}, \eqref{eq:Policy_q}, and parameterized by vector-valued variables $\alpha_p$, $\alpha_q$, respectively.
The control policy parameters $\alpha_p$, $\alpha_q$ are the degrees of freedom of Problem \eqref{eq:AC_sOPF} and can be interpreted as automatic generation control coefficients.
The control policies are generic but not arbitrary because the power balance has to hold.
For the \dc setting (piece-wise) affine control policies can be optimal \cite{Bienstock14,Vrakopoulou17}.
To the best of the authors' knowledge, it is an open research question what policies are optimal for the \ac setting.

Table \ref{tab:Comparison} summarizes the differences between \hopf and \sopf.
In essence, \sopf achieves computational tractability because constraint satisfaction \emph{for every realization} is sacrificed by the introduction of chance constraints.
In general, the optimal solutions of \hopf and \sopf cannot be expected to be equivalent, which leads to a \emph{price of uncertainty}.
The monetary price of uncertainty, for example, may be inferred from the different optimal distributions of costs.
\begin{table}[t]
	\centering
	\caption{Comparison between features of \hopf and \sopf.}
	\label{tab:Comparison}
	\footnotesize
	\renewcommand{\arraystretch}{1.3}
	\begin{tabular}{m{1.65cm}m{2.925cm}m{2.925cm}}
		\toprule
		& \hopf                                                                      & \sopf                                                                                                                                                                \\ \midrule
		DoFs                     & power generation $\pvar{}, \qvar{}$                                        & policy parameters $\alpha_p, \alpha_q$                                                                                                                                    \\
		Optimality               & per-sample optimal                                                         & optimal policies                                                                                                                                                     \\
		Solution characteristics & $\rv{p}^{g \star}$,
		
		$\rv{q}^{g \star}$,
		
		no functional dependency & $\rv{p}^{g \star} = \eta_p(\rvpdvar{}, \rvqdvar{}; \alpha_p^\star)$,
		
		$\rv{q}^{g \star} = \eta_q(\rvpdvar{}, \rvqdvar{}; \alpha_q^\star)$,
		
		functional dependency \\
		Computation              & one run per sample                                                         & single run                                                                                                                                                           \\
		Power balance            & satisfied                                                                  & satisfied                                                                                                                                                            \\
		Constraints              & satisfied per sample                                                       & chance-constraint sense                                                                                                                                              \\ \bottomrule
	\end{tabular}
	\vspace{\adjustlength}
\end{table}
It is worth asking whether the solutions of \hopf and \sopf can be identical?
If so, what are necessary and sufficient conditions?
In this desirable case, a single numerical run provides the solution (cf. \sopf) that is known to strictly satisfy the inequality constraints for all realizations of uncertain generation/demand (cf. \hopf).
The numerical and implementation challenges of online optimization would then be alleviated.

\section{Equivalence of \hopf and \sopf}
We employ polynomial chaos expansion (\pce) as a tool to derive conditions that lead to equivalence between \hopf and \sopf.
Polynomial chaos is a spectral expansion technique for random variables with finite variance that allows to represent a random variable entirely by its deterministic, real-valued \pce coefficients.
With respect to \sopf polynomial chaos has the following advantages:
(i) a stochastic problem is reformulated as a deterministic problem in terms of the deterministic, real-valued \pce coefficients, and
(ii) \pce offers a unified framework for several uncertainty descriptions common in power systems applications, for example Gaussian, Beta, Gamma, and/or Uniform distribution, cf. \cite{Atwa10a, Carpaneto08a, Soubdhan09a}.
The space limitations prohibit a thorough introduction of \pce; we refer to \cite{Sullivan15book,Xiu10book} for further details.
For applications of \pce in the power systems context we refer to \cite{Muehlpfordt17a,Muehlpfordt16b,kit:appino17a,kit:engelmann17b,Ni17, Tang16}.

The next result relies on \pce  to the end of presenting sufficient equivalence conditions for \hopf and \sopf.
\begin{thm}[Equivalence of \hopf and \sopf]~\\
	\label{propo:equivalence}
	Consider the \sopf Problem \eqref{eq:AC_sOPF}, and let the following hold:
	\begin{mylist}
		\item \label{item:DCpowerflow}the \dc power flow assumptions are valid;
		\item \label{item:QuadraticCost} the cost functions are quadratic and positive definite;
		\item \label{item:FinitePCE}the uncertain loads $\rvpdvar{}$, $\rvqdvar{}$ admit a finite \pce that is exact with dimension $L + 1$ with respect to the orthonormal polynomial basis $\{ \basisfun_\ell \}_{\ell = 0}^L$, i.e. $\rvpdvar{} = \sum_{\ell = 0}^{L} \pdvar{\ell} \basisfun_\ell$, where $\pdvar{\ell} \in \mathbb{R}^{\nbus}$;\footnote{This slight abuse of notation w.r.t. \eqref{eq:NetPower_realization} is manageable, because we strictly use the subscript $\ell$ for \pce coefficients in the remainder.}
		\item \label{item:ActiveSet}for all realizations of the uncertain loads $\rvpdvar{}$, $\rvqdvar{}$ the set of active inequality constraints for the respective solution to \eqref{eq:AC_OPF} is the same.
	\end{mylist}
	Then, if the active power policy in the \sopf Problem \eqref{eq:AC_sOPF} is chosen according to
		\begin{align}
		\label{eq:DC_ControlPolicy}
		\rvpvar{} & = \eta_p(\alpha_p) = \sum_{\ell = 0}^{L} \alpha_{p,\ell} \basisfun_\ell,
		\end{align}
	the \sopf solution is identical to the \hopf solution.
	\hfill $\square$
\end{thm}
Before proving the assertions of the theorem, we remark that the seemingly technical assumption \ref{item:FinitePCE} has a clear practical interpretation: for example, any ``canonical'' uncertainty (e.g. Gaussian, Beta, Gamma, Uniform) admits an exact \pce with just two coefficients, i.e. $L + 1 = 2$.
For example, it is sufficient to know the mean and standard deviation of a Gaussian random variable to obtain all statistical moments and its probability density function.
In fact, the \pce of a Gaussian random variable is finite and exact with the mean and standard deviation being the 0\textsuperscript{th} and 1\textsuperscript{st} order \pce coefficient, respectively.\footnote{A basis orthogonal w.r.t. the Gaussian measure on the real line is $\{1, x\}$.}
Note that the active power control policy parameters $\alpha_p$ correspond to the deterministic \pce coefficients of $\rvpvar{}$.
\begin{proof}
Due to space limitations, we present the proof for the case that the set of active inequality constraints from item \ref{item:ActiveSet} is empty.
Under \dc power flow conditions the reactive power is constant due to constant voltage magnitudes, hence not a degree of freedom for the \opf problem.
According to \cite{Muehlpfordt17a} the \sopf Problem \eqref{eq:AC_sOPF}, under assumptions \ref{item:DCpowerflow}--\ref{item:ActiveSet} from Theorem \ref{propo:equivalence}, can be written as a convex quadratic program (\qp)
\begin{subequations}
	\label{eq:DC_sOPF}
	\begin{align}
	\underset{\bs{\alpha}}{\min}\quad  & \frac{1}{2} \, \bs{\alpha}^\top  (I_{L+1} \otimes H) \bs{\alpha} + (e \otimes h)^\top \bs{\alpha}  \\
	\mathrm{s.\, t.} \, ~ ~ ~ \,
	\label{eq:PFE_DC_PCE}
	& (I_{L+1} \otimes \bs{1}^\top) (\bs{\alpha} + \bs{\bs{\pdvar{}}}) = 0,
	\end{align}
\end{subequations}
where $\bs{\alpha} = [\alpha_{p,0}^\top, \hdots, \alpha_{p,L}^\top]^\top$, $\bs{\pdvar{}} = [\pdvar{0}, \hdots, \pdvar{L}]^\top$ are the vectors of stacked \pce coefficients, and $e = [1, 0, \hdots, 0]^\top$ is the $(L+1)$-dimensional unit vector.
The matrix $H \in \mathbb{R}^{\nbus \times \nbus}$ is diagonal with positive entries, hence positive definite.
The \kkt system for \eqref{eq:DC_sOPF} is linear
\begin{equation}
\label{eq:KKT_sOPF}
\bs{A}_{\text{s}} \bs{z}_{\text{s}} = \bs{b}_{\text{s}}.
\end{equation}
The coefficient matrix $\bs{A}_{\text{s}}$, the decision vector $\bs{z}_{\text{s}}$, and right-hand side vector $\bs{b}_{\text{s}}$ correspond to
\begin{align*}
\begin{bmatrix}
(I_{L+1} \otimes H) & (I_{L+1} \otimes \bs{1}^\top)^\top \\
(I_{L+1} \otimes \bs{1}^\top) & 0
\end{bmatrix}
\begin{bmatrix}
\bs{\alpha} \\
\bs{\lambda}
\end{bmatrix}
{=}
- \!
\begin{bmatrix}
(e \otimes h) \\
(I_{L+1} \otimes \bs{1}^\top) \bs{\pdvar{}}
\end{bmatrix}\!,
\end{align*}
where $\bs{\lambda}$ are the \pce coefficients of the Lagrange multiplier for the power balance constraint \eqref{eq:PFE_DC_PCE}, and $I_{L+1}$ is the $(L+1)$ by $(L+1)$ identity matrix.
The coefficient matrix of the linear system \eqref{eq:KKT_sOPF} is regular due to positive definiteness of $H$.
Thus, \eqref{eq:KKT_sOPF} admits a unique solution for $\bs{\alpha}$.
We will now show that \hopf leads to the same system of equations \eqref{eq:KKT_sOPF}, hence leads to the same control policy.
From items \ref{item:DCpowerflow}, \ref{item:QuadraticCost} it follows that the \opf Problem \eqref{eq:AC_OPF} reduces to \dcopf which can be formulated as a convex \qp
\begin{subequations}
	\label{eq:DC_OPF}
	\begin{align}
	\underset{\pvar{}}{\min}\quad  & \frac{1}{2} \, \pvar{}\phantom{}^\top  H \pvar{} + h^\top \pvar{}  \\
	\mathrm{s.\, t.} \, ~ ~ ~ \,
	\label{eq:PFE_DC}
	& \bs{1}^\top (\pvar{} + \pdvar{}) = 0.
	\end{align}
\end{subequations}
The \kkt system for Problem \eqref{eq:DC_OPF} becomes
\begin{equation}
\label{eq:KKT_hOPF}
\begin{bmatrix}
H & \bs{1} \\
\bs{1}^\top & 0 
\end{bmatrix}
\begin{bmatrix}
\pvar{} \\ \lambda
\end{bmatrix}
=
-
\begin{bmatrix}
 h
\\
\bs{1}^\top \pdvar{}
\end{bmatrix},
\end{equation}
where $\lambda \in \mathbb{R}$ is the Lagrange multiplier for the \dc power flow balance \eqref{eq:PFE_DC}.
The coefficient matrix in \eqref{eq:KKT_hOPF} is regular because $H$ is positive definite.
Hence, the argmin operator is linear in the demand $\pdvar{}$, and can be used directly for uncertainty propagation.
To that end, the \pce for the uncertain demand $\rvpdvar{}$ and generation $\rvpvar{}$, see Assumption \ref{item:FinitePCE} and \eqref{eq:DC_ControlPolicy}, are substituted in \eqref{eq:KKT_hOPF}, and the \pce for the Lagrange multiplier $\rv{\lambda} = \sum_{\ell = 0}^{L} \lambda_\ell \basisfun_\ell$ is introduced.
The Galerkin-projected system in matrix form becomes
\begin{subequations}
	\label{eq:KKT_hOPF_Galerkin}
	\begin{align}
	\begin{bmatrix}
	I_{L+1} \otimes \begin{bmatrix}
	H & \bs{1} \\
	\bs{1}^\top & 0 
	\end{bmatrix}
	\end{bmatrix}
	\bs{z}_{\text{h}}
	=
	\bs{A}_{\text{h}} \bs{z}_{\text{h}} = 
	\bs{b}_{\text{h}},
	\end{align}
where $\bs{z}$ contains all \pce coefficients $\alpha_\ell$ and $\lambda_\ell$, and $\bs{b}$ contains the \pce coefficients $\pdvar{\ell}$ and the vector of linear cost coefficients~$h$,
	\begin{align}
	\bs{z}_{\text{h}} &= \begin{bmatrix}
	\alpha_{p,0}\phantom{}^\top &
	\lambda_0 &
	\alpha_{p,1} & 
	\lambda_1 &
	\hdots &
	\alpha_{p,L}\phantom{}^\top &
	\lambda_L
	\end{bmatrix}^\top, \\
	\bs{b}_{\text{h}} &= 
	-
	\begin{bmatrix}
	h^\top & \bs{1}^\top \pdvar{0} & 0 & \bs{1}^\top \pdvar{1} & \hdots & 0 & \bs{1}^\top \pdvar{L}
	\end{bmatrix}^\top.
	\end{align}
\end{subequations}
To show equivalence between \sopf and \hopf from Theorem~\ref{propo:equivalence} it remains to show that the linear systems \eqref{eq:KKT_sOPF} and \eqref{eq:KKT_hOPF_Galerkin} admit the same solution.
To do so, introduce the following permutation matrix
\begin{equation}
M = \begin{bmatrix}
I_{L+1} \otimes \begin{bmatrix}
I_{\nbus} & 0
\end{bmatrix}
\\
I_{L+1} \otimes \begin{bmatrix}
0 & 1
\end{bmatrix}
\end{bmatrix}
\in \mathbb{R}^{(L+1)(n+1) \times (L+1)(n+1) },
\end{equation}
and observe that
\begin{align}
\bs{z}_{\text{s}}
&=
M
\bs{z}_{\text{h}}, ~
\bs{b}_{\text{s}}
= M \bs{b}_{\text{h}}, ~
\bs{A}_{\text{s}} = M^\top \bs{A}_{\text{h}} M.
\end{align}
In other words, the linear systems \eqref{eq:KKT_sOPF} and \eqref{eq:KKT_hOPF_Galerkin} are equivalent after permuting with $M$.

Recall that in \eqref{eq:AC_OPF} we consider box constraints. Moreover, observe that whenever the set of active inequality constraints is nonempty, yet does not change for all realizations, the chance constraints in  \eqref{eq:AC_sOPF} are not active. Instead they are replaced by active linear inequalities, i.e. by linear equalities. Thus, the same steps as above prove the assertion.
\end{proof}
We remark that the optimal solution to the linear systems \eqref{eq:KKT_sOPF}, \eqref{eq:KKT_hOPF_Galerkin} is affine in the \pce coefficients of the demand.
If the active set changes for some realizations, the optimal solution can still be parameterized affinely for each active set.

Loosely speaking, Theorem~\ref{propo:equivalence} states conditions that lead to the price of uncertainty being zero.
That in turn means that if the conditions from Theorem~\ref{propo:equivalence} are satisfied, no online optimization scheme will yield better solutions than the \ccopf solution.
It is worth asking how the solutions from \sopf and \hopf can be compared in case the assumptions of Theorem~\ref{propo:equivalence} do not hold, i.e. in case of a non-zero price of uncertainty.
Importantly, the case when constraints are active for just a few number of realizations $\pdvar{}$ of the uncertain demand $\rvpdvar{}$.
Recall that the solutions stemming from \sopf and \hopf are random variables.
Hence, their probability density functions (\pdfs) can be compared using different metrics.
In statistics, popular choices include the Kullback-Leibler divergence or the Hellinger distance \cite{Gibbs02}.
However, the Kullback-Leibler divergence may not be meaningful, for instance, in case of compact supports that overlap.
The Hellinger distance leads to numerical difficulties when a distribution contains Dirac-deltas---which can be the case for \hopf when constraints are active.
In the present paper, we suggest relying on the total variational distance \cite{Gibbs02}.
For two random variables $\rv{x}$ and $\rv{y}$ with \pdfs $f_{\rv{x}}$ and $f_{\rv{y}}$, the total variational distance is given by $\Delta: \mathrm{L}^2(\Omega_x, \mu_x; \mathbb{R}) \times \mathrm{L}^2(\Omega_y, \mu_y; \mathbb{R}) \to [0,1]$ 
\begin{equation}
\Delta(\rv{x}, \rv{y}) = \frac{1}{2} \, \int_{\mathbb{R}} | f_{\rv{x}}(\tau) - f_{\rv{y}}(\tau)| \mathrm{d}\tau.
\end{equation}
The total variational distance $\Delta$ can serve as an indicator, for example, as when to prefer \ccopf to fast online optimization.
A small value of $\Delta$ indicates that \ccopf and \hopf lead to similar power injections.
Thus, computing the \ccopf once and applying it online, will yield good results.
In contrast, a large value of $\Delta$ indicates that fast online optimization approaches will outperform the \ccopf solution in terms of constraint satisfaction and cost.

\section{Tutorial Example}
\begin{figure}
	\centering
	\includegraphics[]{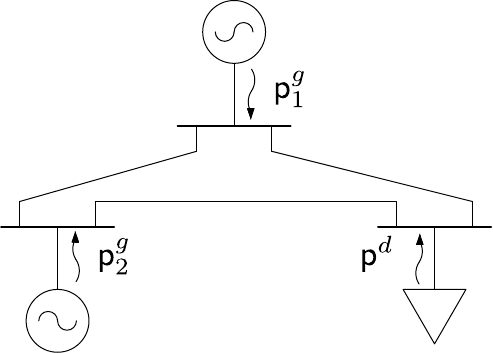}
	\caption{Three-bus network with two generators and one load.}
	\label{fig:ThreeBus}	
	\vspace{\adjustlength}
\end{figure}
The following tutorial example demonstrates Theorem~\ref{propo:equivalence} in action, and quantifies the price of uncertainty in case of \hopf and \sopf yielding different solutions.
We focus on the rather small three-bus example because it is simple enough to have an analytical solution, and it is complex enough to provide helpful insight.
All units are in per-unit values.

Consider the connected three-bus network from Figure \ref{fig:ThreeBus}, which has a generator but no load connected to buses~1 and~2, respectively, and a load but no generator connected to bus~3.
With slight abuse of notation we set $\rvpdvar{} \equiv \rvpdvar{3}$, and $\rvpvar{} \equiv [\rvpvar{1}, \rvpvar{2}]^\top$.
The considered deterministic \dcopf is
\begin{subequations}
	\label{eq:DCOPF_Example}
	\begin{align}
	\underset{\pvar{}}{\min}\quad  & \frac{1}{2} \, \pvar{}\phantom{}^\top  \begin{bmatrix}
		H_{11} & 0 \\
		0 & H_{22}
	\end{bmatrix} \pvar{} + \begin{bmatrix}
	h_1 & h_2
	\end{bmatrix}^\top \pvar{}  \\
	\mathrm{s.\, t.} \, ~ ~ ~ \,
	& \pvar{1} + \pvar{2} + \pdvar{} = 0, \\
	& \pvar{1} \leq p_1^{\text{max}},
	\end{align}
\end{subequations}
with $H_{11}, H_{22} > 0$.
Note that power demand is counted negative.
The argmin operator to Problem \eqref{eq:DCOPF_Example} is 
\begin{equation}
\label{eq:DCOPF_example_argmin}
\mathbb{R}^2 \ni \pvar{}\phantom{}^\star = 
\begin{cases}
\begin{bmatrix}
-1 \\ \phantom{-}1
\end{bmatrix} \beta
-
\begin{bmatrix}
\gamma \\
1 - \gamma
\end{bmatrix}
\pdvar{}, &- \pdvar{} < \frac{p_1^{\text{max}} + \beta}{\gamma}, \\
\begin{bmatrix}
p_1^{\text{max}} \\
-(\pdvar{} + p_1^{\text{max}})\end{bmatrix},
& -\pdvar{} \geq \frac{p_1^{\text{max}} + \beta}{\gamma},
\end{cases}
\end{equation}
with $\beta {=} (h_1 {-} h_2)/(H_{11} {+} H_{22})$, and $\gamma {=} H_{22}/(H_{11} + H_{22}) {>} 0$.

If the probability density function $f_{\rvpdvar{}}$ of the uncertain demand $\rvpdvar{}$ is given, the density of the \hopf solution is obtained from the argmin~\eqref{eq:DCOPF_example_argmin} as
\begin{align}
\nonumber
&f_{\rvpvar{1}}(x_1) {=} \begin{cases}
\frac{1}{\gamma} \, f_{\rvpdvar{}}\!\left(\frac{x_1 + \beta}{-\gamma}\right), & x_1 < p_1^{\text{max}}, \\
\left(1 - \frac{1}{\gamma} F_{\rvpdvar{}}\!\! \left( \frac{p_1^{\text{max}} + \beta}{-\gamma} \right)\right) h(p_1^{\text{max}} {-} x_1), & x_1 = p_1^{\text{max}},
\end{cases} \\
\label{eq:PDF_HOPF}
&f_{\rvpvar{2}}(x_2) {=} \begin{cases}
\frac{1}{1 - \gamma} \, f_{\rvpdvar{}}\!\!\left( \frac{x_2 - \beta}{\gamma - 1} \right), & x_2 >   \frac{p_1^{\text{max}} + \beta}{\gamma} - p_1^{\text{max}}, \\
f_{\rvpdvar{}}(-x_2 - p_1^{\text{max}}), & x_2 \leq \frac{p_1^{\text{max}} + \beta}{\gamma} - p_1^{\text{max}},
\end{cases}
\end{align}
where $h$ is the Dirac-delta.
The case-dependent definition of the \pdfs is due to the upper generation limit $p_1^{\text{max}}$.
If generator~1 operates below its limit $p_1^{\text{max}}$, the share of active power generation assigned to generator~1 and~2 is determined by the cost coefficients $H_{11}, H_{22}, h_1, h_2$ via $\beta$ and $\gamma$; this situation corresponds to the upper cases in \eqref{eq:PDF_HOPF}.
At the generation limit $p_1^{\text{max}}$ the \pdf of generator 1 becomes a delta-pulse with ``height'' equal to the mass of the \pdf that is cut off to the right of $p_1^{\text{max}}$.
In case of generator~1 hitting its limit, generator 2 supplies the remaining active power to meet the power demand, and to guarantee power balance.
To this end, the \pdf of power generation at bus 2 has a discontinuity and becomes equivalent to the shifted \pdf of the demand, $f_{\rvpdvar{}}(-x_2 - p_1^{\text{max}})$.

\newpage
We now turn to the solution via \sopf.
We simplify the notation of the policy paramter to $\alpha \leftarrow \alpha_p$.
In light of Theorem~\ref{propo:equivalence}, Assumption~\ref{item:FinitePCE}, we assume that the \pce for the demand $\rvpdvar{}$ is finite and exact with $L + 1 = 2$, and respective \pce coefficients $\pdvar{0}$, $\pdvar{1}$.
For example, this is the case for any ``canonical'' uncertainty in the corresponding basis, i.e. Gaussian, Beta, Gamma, or Uniform distribution.
According to Theorem~\ref{propo:equivalence} the active power policy will consist of two elements, too.
Let $\alpha^\star_{\ell}$ denote the optimal solution for the $\ell$\textsuperscript{th} coefficient to Problem~\eqref{eq:AC_sOPF} tailored to the setting from Problem~\eqref{eq:DCOPF_Example}.
Then, the \pdf of the optimal active power policies for generators~1 and 2 is
\begin{equation}
\label{eq:PDF_PCE}
f_{\rvpvar{i}}(\pvar{i}) = \left|\frac{\pdvar{1}}{\alpha_{1,i}^\star} \right|  \, f_{\rvpdvar{}}\!\left( \frac{\alpha^\star_{1,i} \pdvar{0} - \alpha^\star_{0,i} \pdvar{1}}{\alpha^\star_{1,i}} + \frac{\pdvar{1} }{\alpha^\star_{1,i}}\, \pvar{i} \right),
\end{equation}
for $i = 1,2$, where $\alpha^\star_{\ell,i}$ denotes the $i$\textsuperscript{th} entry of the $\ell$\textsuperscript{th} coefficient $\alpha^\star_{\ell}$.
Compared to the \pdf from \hopf, the \pdf \eqref{eq:PDF_PCE} from \sopf is always continuous, possibly leading to a price of uncertainty, because the generation limit may be violated (with user-specified low probability).
\begin{table}[t]
	\centering
	\caption{Numerical values for Problem~\eqref{eq:DCOPF_Example}.	\label{tab:Numericalvalues}}
	\footnotesize
	\begin{tabular}{cccccc}
		\toprule
		Case   &    $H_{11}$    & $H_{22}$ & $h_1$ & $h_2$ & $p_1^{\text{max}}$ \\ \midrule
		\caseunc & $\{0.2, 0.3\}$ &  $0.2$   & $0.5$ & $0.6$ &       $1.5$        \\
		\casecon &     $0.2$      &  $0.2$   & $0.5$ & $0.6$ &       $0.85$       \\ \bottomrule
	\end{tabular}
	\vspace{\adjustlength}
\end{table}
Assume now, in light of item~\ref{item:ActiveSet} of Theorem~\ref{propo:equivalence}, that for all realizations $\pdvar{}$ of the uncertain demand $\rvpdvar{}$ it holds that $- \pdvar{} < (p_1^{\text{max}} + \beta)/\gamma$.
In other words, the demand never results in the generation constraint being active.
The \pce coefficients of the optimal active power generation are obtained from the argmin operator \eqref{eq:DCOPF_example_argmin}
\begin{equation}
\label{eq:Example_PCE_Solution}
\mathbb{R}^2 \ni
\alpha_{\ell}^\star =
\begin{cases}
\begin{bmatrix}
-1 \\ \phantom{-}1
\end{bmatrix} \beta
-
\begin{bmatrix}
\gamma \\
1 - \gamma
\end{bmatrix}
\pdvar{0}, & \ell = 0, \\
\phantom{\begin{bmatrix}
	-1 \\
\phantom{-}1
	\end{bmatrix} \beta}
- \begin{bmatrix}
\gamma \\
1 - \gamma
\end{bmatrix}
\pdvar{\ell}, & \ell = 1, \hdots, L.
\end{cases}
\end{equation}
The coefficients from \eqref{eq:Example_PCE_Solution} are required to determine the \pdf from \sopf according to \eqref{eq:PDF_PCE}.

\begin{figure*}
	\centering
	\subfloat[\pdf of demand $\rvpdvar{}$.\label{fig:demand}]{%
		\includegraphics[width=3.5cm]{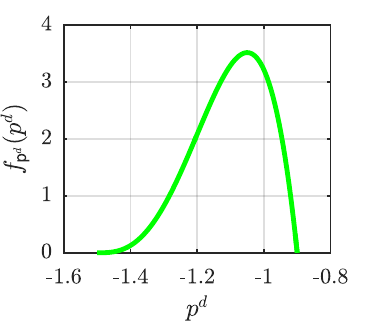}%
	}
	\subfloat[Case \caseunc yields zero price of uncertainty.\label{fig:equivalence}]{
		\includegraphics[width=3.5cm]{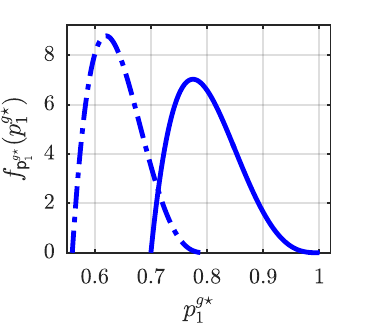}%
		\includegraphics[width=3.5cm]{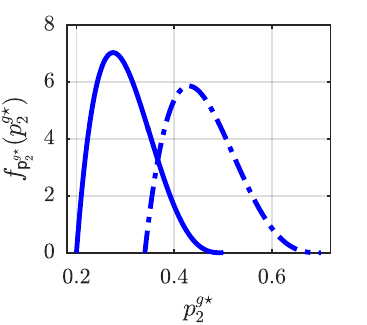}%
	}	
	\subfloat[Case \casecon yields non-zero price of uncertainty.\label{fig:nonequivalence}]{%
		\includegraphics[width=3.5cm]{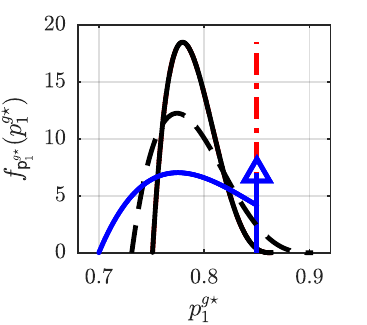}%
		\includegraphics[width=3.5cm]{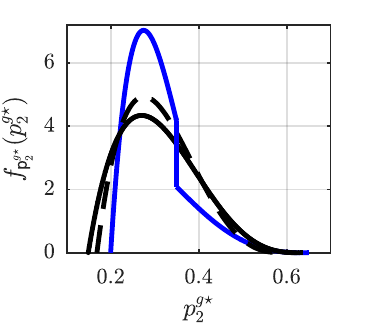}%
	}	
	\vspace{2mm}			
	\caption{Probability density functions of demand and generation for cases given in Table~\ref{tab:Numericalvalues}. Blue denotes \hopf, black denotes \sopf.}
	\vspace{\adjustlength}
\end{figure*}

Specifically, let the uncertain demand follow a Beta distribution with support $[p^{d,\text{min}}, p^{d,\text{max}}] = [-1.5, -0.9]$, and shape parameters $a = 4$, $b = 2$, i.e. $\rvpdvar{} \sim \mathrm{B}([-1.5, -0.9], 4, 2)$. The skewed \pdf is shown in Figure~\ref{fig:demand}.
In the respective Jacobi polynomial basis the \pce for $\rvpdvar{}$ is finite and exact with the following two \pce coefficients
$ \pdvar{0} = -1.1$, $\pdvar{1} = 0.1$.
In case of all assumptions of Theorem~\ref{propo:equivalence} being satisfied, the densities from \hopf, \eqref{eq:PDF_HOPF}, and \sopf, \eqref{eq:PDF_PCE}, are equivalent.
Consider case \caseunc from Table \ref{tab:Numericalvalues}, for which all realizations satisfy 
\begin{equation}
\label{eq:ConditionForEquivalence}
- \pdvar{}  \leq 1.5 < (p_1^{\text{max}} + \beta)/\gamma \in  \{2.50, 3.25 \}.
\end{equation}
Figure~\ref{fig:equivalence} shows the optimal densities for active power generation.
The share of generation is entirely determined by the cost coefficients from Table~\ref{tab:Numericalvalues} that enter the optimal solution~\eqref{eq:PDF_HOPF} via $\beta$ and $\gamma$.
The case $H_{11} = 0.2$, shown in solid blue in Figure~\ref{fig:equivalence}, leads to significantly higher generation at bus~1 compared to the case $H_{11} = 0.3$, shown dash-dotted blue in Figure~\ref{fig:equivalence}.
The active power limit $p_1^{\text{max}} = 1.5$ is not attained.
Equivalence of the optimal solutions means that the price of uncertainty is zero.
\begin{table}[t]
	\centering
	\caption{Constraint satisfaction and total variational distance.\label{tab:SimulationParameters}}
	\footnotesize	
	\begin{tabular}{cccc}
		\toprule
		$\delta$ & $\pr{\rvpvar{1} \leq p_1^{\text{max}}}$ & $\Delta (\rvpvar{1,\text{\hopf}}, \rvpvar{1,\text{\sopf}})$ & $\Delta (\rvpvar{2,\text{\hopf}}, \rvpvar{2,\text{\sopf}})$ \\ \midrule
		   2.0     &                96.51\,\%                &                        0.3197                         &                        0.1882                         \\
		   3.0     &                99.86\,\%                &                        0.4734                         &                        0.2451                         \\ \bottomrule
	\end{tabular}
	\vspace{\adjustlength}
\end{table}

Consider now case \casecon from Table~\ref{tab:Numericalvalues}, instead, for which condition \eqref{eq:ConditionForEquivalence} is violated
$
- \pdvar{}  \leq 1.5 \not< (p_1^{\text{max}} + \beta)/\gamma = 1.2.
$
The \pdf~\eqref{eq:PDF_HOPF} from \hopf is plotted in blue in Figure~\ref{fig:nonequivalence}.
It shows the discontinuity at the upper limit $p_1^{\text{max}} = 0.85$ where the \pdf becomes a delta-pulse---denoted by the triangle in Figure~\ref{fig:nonequivalence}.
Effectively, the \pdf for bus~1 from the unconstrained case~\caseunc is cut off at the generation limit, and bus~2 accounts for the remainder.
To obtain the solution via \sopf the chance constraint $\pr{\rvpvar{1} \leq p_1^{\text{max}}}$ is reformulated as $\ev{\rvpvar{1}} + \delta \, \Sigma\left[\rvpvar{1}\right] \leq p_1^{\text{max}}$. The expected value $\ev{\rvpvar{1}}$ and the standard deviation $\Sigma[\cdot]$ can be computed from the \pce coefficients \cite{Sullivan15book}.
The continuous \pdf \eqref{eq:PDF_PCE} from \sopf violates the constraint limit, in contrast to the \pdf via \hopf; in Figure~\ref{fig:nonequivalence} the dashed black plot is for $\delta = 2$, and the solid black plot is for $\delta = 3$.
Qualitatively, a higher value of $\delta$ leads the \sopf solution to stay away from the constraint.
Quantitatively, the probabilities for constraint satisfaction are summarized in Table \ref{tab:SimulationParameters}.
The parameter $\delta$ has another influence on the quality of the solution:
Compared to the unconstrained case~\caseunc, the \pdf for bus~1 in case~\casecon is significantly more narrow and the mode is shifted to values of higher injection; the opposite effect can be observed for the power generation at bus~2 which becomes more wide.
This leads to less variability in the generation at bus~1, which necessarily leads to higher variability in power generation at bus~2 to ensure power balance.
The higher the value for $\delta$, the less variability is allowed at bus~1.

The value of $\delta$ also affects the total variational distance $\Delta(\cdot, \cdot)$; Table~\ref{tab:SimulationParameters} lists the numerival values.
The more narrow \pdf at bus~1 for $\delta = 3$ leads to a 48\,\% increase in the total variational distance compared to $\delta = 2$.
A similar behavior is observed at bus~2, where the ``true'' \pdf from \hopf is fairly narrow, but the \pdfs from \sopf are structurally too wide.
In that case, for $\delta = 3$ the total variational distance is 30\,\% larger compared to $\delta = 2$.

\section{Conclusion \& Outlook}
This paper relates chance-constrained \opf to  in-hindsight \opf, which serves as a full-information yet unrealistic benchmark.
For \ccopf, an entire control policy is computed by means of a single optimization problem \emph{before} the realization of the uncertainty is known; at the expense of possible constraint violations (with user-specified low probability).
For \hopf, an \opf problem is solved for every realization of the uncertainty, which leads to per-sample optimality, infinitely many \opf problems, and no immediate control policy by means of a functional dependency.
We show that \ccopf and \hopf are equivalent for \dcopf problems for which the active set of inequality constraints is unchanged for all realizations of the uncertainty.
In that case, the policy from \ccopf gives equivalent results to online optimization approaches.
If the solutions from \ccopf and \hopf do differ, the size of the total variational distance may indicate whether \ccopf should be favored to fast online optimization, or vice versa.
A tutorial three-bus example underpins our results.

The present paper discusses elements of the relation of \ccopf, \hopf, and fast online optimization.
However, several questions remain open: What other dimensions enter the price of uncertainty? Can a measurement-based  detection of changes in the active set trigger new \ccopf computations? What can be said in the multi-stage \ac setting? 

\printbibliography
\end{document}